\documentclass[11pt]{article}
\usepackage{amssymb}
\usepackage{amsmath}
\usepackage{t1enc}
\usepackage[latin1]{inputenc}
\usepackage[english]{babel}
\usepackage{amsmath,amsthm}
\usepackage{amsfonts}
\usepackage{latexsym}
\usepackage[dvips]{graphicx}
\usepackage{graphicx}
\usepackage[natural]{xcolor}
\usepackage{setspace}
\usepackage{fullpage}

\newtheorem{definition}{Definition}
\newtheorem{theorem}{Theorem}
\newtheorem{lemma}{Lemma}
\newtheorem{proposition}{Proposition}

\newtheorem{remark}{Remark}
\newtheorem{corollary}{Corollary}
\makeatletter
   
   \@addtoreset{equation}{section}
\makeatother

\begin{document}

\author{Liangquan Zhang$^{1,2}\thanks{
This work was partially supported by the National Nature Science Foundation
of China (No. 11201263) and the Nature Science Foundation of Shandong
Province (No. ZR2012AQ004) and the Natural Science Foundation of China (No.
11201264) and the ITEA MODRIO project of INRIA. E-mail: xiaoquan51011@163.com.}$ \\
{\small 1. School of Science, }\\
{\small Beijing University of Posts and Telecommunications, }\\
{\small Beijing 100876, China}\\
{\small 2. INRIA, Campus de Beaulieu, 35042 Rennes Cedex, France }}
\title{The Viability Property for Path-dependent SDE under Open Constraints}
\maketitle

\begin{abstract}
In this note, we study the viability of a bounded open domain in $\mathbb{R}%
^{n}$ for a process driven by a path-dependent stochastic differential
equation with Lipschitz data. We extend an invariant result of Cannarsa, Da.
Prato and Frankowska [\textit{Indiana Univ. Math. J.} \textbf{59} (2010)
53-78] to a non-Markovian setting.
\end{abstract}



\section{Introduction}

Let $v=\left( \Omega ,\mathcal{F},\left( \mathcal{F}_{t}\right) _{t\geq
0},P,W\right) $ be a reference probability system composed of a completed
probability space $\left( \Omega ,\mathcal{F},P\right) ,$ a filtration $%
\left( \mathcal{F}_{t}\right) _{t\geq 0}$ satisfying the usual assumptions
of right-continuity and completeness, and a $d$-dimensional $\left( \mathcal{%
F}_{t}\right) $-Brownian motion $W$ defined on $\left( \Omega ,\mathcal{F}%
,P\right) .$

Let $U$ denote an open domain in $\mathbb{R}^{n}.$ In particular, consider a
continuous martingale $X$ governed by a stochastic differential equation
(SDE in short) in the sense that the coefficients of the SDE are allowed to
depend on the path of the process which, of course, is in a non-Markovian
setting. More precisely, the SDE is defined as follows:
\begin{equation}
\left\{
\begin{array}{rcl}
\mbox{\rm d}X^{\chi _{t}}\left( r\right) & = & \mu \left( X_{r}^{\chi
_{t}}\right) \mbox{\rm d}r+\sigma \left( X_{r}^{\chi _{t}}\right)
\mbox{\rm
d}W\left( r\right) ,\text{\qquad }t<r, \\
\text{ }X^{\chi _{t}}\left( s\right) & = & \chi _{t}\left( s\right) ,\text{%
\qquad if }0\leq s\leq t,%
\end{array}%
\right.  \label{1.1}
\end{equation}%
where $\chi _{t}\in \Lambda _{U_{t}}$ which denotes all the right-continuous
functions with left limits defined on $\left[ 0,t\right] $ valued in $U$. $%
\mu $ and $\sigma $ are functionals on $\Lambda $ which denotes all the
right-continuous functions with left limits defined on $\left[ 0,t\right] $
for $\forall t>0$, valued in $\mathbb{R}^{n}$. Besides, denote by $\chi
\left( t\right) $ the value of $\chi $ at $t$ and by $\chi _{t}$ the
restriction on $\left[ 0,t\right] ,$ the same for $X^{\chi _{t}}.$

A natural question arises$,$ if $\chi _{t}\in \Lambda _{U_{t}},$ under what
assumptions can one claim
\begin{equation*}
X^{\chi _{t}}\left( r\right) \in U,\qquad r\in \left[ t,+\infty \right)
\text{ }P\text{-a.s. }?
\end{equation*}%
We study this problem by means of viability theory (see \cite{A, AF} for
more details).

Let $\mathcal{K}$ be a closed subset of $\mathbb{R}^{n}$. We say that $%
\mathcal{K}$ is viable for (\ref{1.1}) if and only if for any $\chi _{t}\in
\Lambda _{\mathcal{K}_{t}}$, the solution to (\ref{1.1}) satisfies
\begin{equation}
X^{\chi _{t}}\left( s\right) \in \mathcal{K},\quad \forall s\geq t,\text{ }P%
\text{-almost surely.}  \label{1.2}
\end{equation}%
In fact, the property of viability of $\mathcal{K}$ for the classical
systems have been extensively studied. We refer the reader to the monographs
\cite{A} for the deterministic case. For the stochastic case, several
results have been obtained: through stochastic tangent cones in \cite{A}-%
\cite{AD}, through viscosity solutions of partial differential equations in
\cite{BCQ}-\cite{BJ} (for more information see references therein). We
mention that in \cite{CPF}, in order to study the existence and uniqueness
of the invariant measure associated to the transition semigroup of a
diffusion process in a bounded open subset of $\mathbb{R}^{n}$, the authors
investigated the invariance of a bounded open domain with piecewise smooth
boundary for a Markov system. Motivated by this paper, we consider the
viability for (\ref{1.1}) under some open smooth domain.

Recently, Dupire \cite{D} 2009 defined a notion of directional derivative
for functionals. Using this notion, Cont and Fournié \cite{CF2} extended Fö%
llmer's (see \cite{F1}) pathwise change of variable formula to
non-anticipative functionals on the space of cadlag paths. Their results
lead to functional extension of the Itô's formula for continuous
semimartingale on open set (see Proposition 7, in \cite{CF2}). Based on
above work, we may use the extension of the Itô's formula to study the
viability in the framework of path-dependence systems.

The rest of this paper is organized as follows: after some preliminaries in
the second section, we are devoted the third section to developing the
viability result for some smooth domain.

\section{Preliminaries and Notations}

Throughout this paper, the notations are mainly taken from Cont and Fournié
\cite{CF1, CF2, CF3} and Dupire \cite{D}. For a cadlag path $x\in D\left( %
\left[ 0,T\right] ,\mathbb{R}^{n}\right) ,$ denote by $x\left( t\right) $
the value of $x$ at $t$ and by $x_{t}=\left( x\left( u\right) ,0\leq u\leq
t\right) $ the restriction of $x$ to $\left[ 0,t\right] .$ Thus $x_{t}\in
D\left( \left[ 0,t\right] ,\mathbb{R}^{n}\right) .$ Similarly, for a
stochastic process $X$ we shall denote $X\left( t\right) $ its value at $t$
and $X_{t}=\left( X\left( u\right) ,0\leq u\leq t\right) $ its path on $%
\left[ 0,t\right] .$

Let $T>0$ be a fixed time horizon and $U\subset \mathbb{R}^{n}$ be an open
subset of $\mathbb{R}^{n}$ and $S\subset \mathbb{R}^{m}$ be a Borel subset
of $\mathbb{R}^{m}$. We denote the boundary of $U$ by $\partial U$, the
closure of $U$ by $\overline{U}=U\cup \partial U.$ We call "$U$-valued
cadlag function"\ a right continuous $f:\left[ 0,T\right] \mapsto U$ with
left limits such that for each $t\in \left[ 0,T\right] ,$ $f\left( t-\right)
\in U.$ Denote by $\mathcal{U}_{t}=D\left( \left[ 0,t\right] ,U\right) $
(resp. $\mathcal{S}_{t}=D\left( \left[ 0,t\right] ,S\right) $) the space of $%
U$-valued cadlag functions (resp. $S$), and $C_{0}\left( \left[ 0,t\right]
,U\right) $ the set of continuous functions with values in $U.$

\begin{definition}[Non-anticipative functionals on path space]
A non-anticipative functional on $\mathcal{U}_{T}$ is a family $F=\left(
F_{t}\right) _{t\in \left[ 0,T\right] }$ of maps%
\begin{equation*}
F_{t}:\mathcal{U}_{t}\rightarrow \mathbb{R}.
\end{equation*}
\end{definition}

We consider throughout this paper non-anticipative functionals%
\begin{equation*}
F=\left( F_{t}\right) _{t\in \left[ 0,T\right] },\qquad F_{t}:\mathcal{U}%
_{t}\times \mathcal{S}_{t}\rightarrow \mathbb{R},
\end{equation*}%
where $F$ has a "predictable"\ dependence with respect to the second
argument:%
\begin{equation}
\forall t\in \left[ 0,T\right] ,\text{ }\forall \left( x,v\right) \in
\mathcal{U}_{t}\times \mathcal{S}_{t},\qquad F_{t}\left( x_{t},v_{t}\right)
=F_{t}\left( x_{t},v_{t-}\right) .  \label{2.1}
\end{equation}%
$F$ can be view as a functional on the vector bundle $\Upsilon ^{U\times
S}=\cup _{t\in \left[ 0,T\right] }\mathcal{U}_{t}\times \mathcal{S}_{t}.$

For each $\gamma _{t}\in \mathcal{U}_{t}$ denote%
\begin{equation*}
\gamma _{t,\delta }\left( s\right) =\left\{
\begin{array}{ll}
\gamma _{t}\left( s\right) , & \text{for }0\leq s<t, \\
\gamma _{t}\left( t\right) , & \text{for }t\leq s\leq t+\delta ,%
\end{array}%
\right.
\end{equation*}%
and for $x\in \mathbb{R}^{n}$ small enough,
\begin{equation*}
\gamma _{t}^{x}\left( s\right) =\left\{
\begin{array}{ll}
\gamma \left( s\right) , & \text{for }0\leq s<t, \\
\gamma \left( t\right) +x, & \text{for }s=t.%
\end{array}%
\right.
\end{equation*}%
Apparently, $\gamma _{t,\delta }\in \mathcal{U}_{t+\delta }.$ We also denote
\begin{equation*}
\left( \gamma _{t}^{x}\right) _{t,\delta }\left( s\right) =\left\{
\begin{array}{ll}
\gamma \left( s\right) , & \text{for }0\leq s<t, \\
\gamma \left( t\right) +x, & \text{for }t\leq s\leq t+\delta .%
\end{array}%
\right.
\end{equation*}%
Let $\bar{\gamma}_{\bar{t}},$ $\gamma _{t}\in \mathcal{U}_{t}$ be given with
$\bar{t}\leq t,$ we denote $\bar{\gamma}_{\bar{t}}\oplus \gamma _{t}\in
\mathcal{U}_{t}$ by
\begin{equation*}
\bar{\gamma}_{\bar{t}}\oplus \gamma _{t}=\left\{
\begin{array}{ll}
\bar{\gamma}_{\bar{t}}\left( s\right) , & \text{for }0\leq s<\bar{t}, \\
\gamma _{t}\left( s\right) , & \text{for }\bar{t}\leq s<t.%
\end{array}%
\right.
\end{equation*}%
We now introduce a distance between two paths, not necessarily defined on
the same time interval. For $T\geq t^{\prime }=t+h\geq t\geq 0,$ $\left(
x,v\right) \in \mathcal{U}_{t}\times \mathcal{S}_{t}$ and $\left( x^{\prime
},v^{\prime }\right) \in D\left( \left[ 0,t+h\right] ,\mathbb{R}^{n}\right)
\times \mathcal{S}_{t+h}$ define%
\begin{eqnarray*}
\left\Vert x\right\Vert &=&\sup\limits_{0\leq s\leq t}\left\vert x_{t}\left(
s\right) \right\vert , \\
d_{\infty }\left( \left( x,v\right) ,\left( x^{\prime },v^{\prime }\right)
\right) &=&\sup\limits_{r\in \left[ 0,t+h\right] }\left\vert x_{t,h}\left(
r\right) -x^{\prime }\left( r\right) \right\vert +\sup\limits_{r\in \left[
0,t+h\right] }\left\vert v_{t,h}\left( r\right) -v^{\prime }\left( r\right)
\right\vert +h.
\end{eqnarray*}%
$\left( \Upsilon ^{U\times S},d_{\infty }\right) $ is a metric space.

\begin{definition}[Continuous at fixed times]
A non-anticipative functional $F=\left( F_{t}\right) _{t\in \left[ 0,T\right]
}$ is said to be continuous at fixed times if for any $t\leq T,$ $F_{t}:%
\mathcal{U}_{t}\times \mathcal{S}_{t}\rightarrow \mathbb{R}$ is continuous
for the supremum norm.
\end{definition}

\begin{definition}[Left-continuous functionals]
Define $\mathbb{F}_{l}^{\infty }$ as the set of functionals $F=\left(
F_{t},t\in \left[ 0,T\right] \right) $ which verify:%
\begin{equation*}
\begin{array}{c}
\forall t\in \left[ 0,T\right] ,\text{ }\varepsilon >0,\text{ }\forall
\left( x,v\right) \in \mathcal{U}_{t}\times \mathcal{S}_{t},\text{ }\exists
\eta >0,\text{ }\forall h\in \left[ 0,t\right] , \\
\forall \left( x^{\prime },v^{\prime }\right) \in \mathcal{U}_{t-h}\times
\mathcal{S}_{t-h},\text{ }d_{\infty }\left( \left( x,v\right) ,\left(
x^{\prime },v^{\prime }\right) \right) <\eta \Rightarrow \left\vert
F_{t}\left( x,v\right) -F_{t-h}\left( x^{\prime },v^{\prime }\right)
\right\vert <\varepsilon .%
\end{array}%
\end{equation*}
\end{definition}

\begin{definition}[Right-continuous functionals]
Define $\mathbb{F}_{r}^{\infty }$ as the set of functionals $F=\left(
F_{t},t\in \left[ 0,T\right] \right) $ which verify:%
\begin{equation*}
\begin{array}{c}
\forall t\in \left[ 0,T\right] ,\text{ }\varepsilon >0,\text{ }\forall
\left( x,v\right) \in \mathcal{U}_{t}\times \mathcal{S}_{t},\text{ }\exists
\eta >0,\text{ }\forall h\in \left[ 0,T-t\right] , \\
\forall \left( x^{\prime },v^{\prime }\right) \in \mathcal{U}_{t+h}\times
\mathcal{S}_{t+h},\text{ }d_{\infty }\left( \left( x,v\right) ,\text{ }%
\left( x^{\prime },v^{\prime }\right) \right) <\eta \Rightarrow \left\vert
F_{t}\left( x,v\right) -F_{t+h}\left( x^{\prime },v^{\prime }\right)
\right\vert <\varepsilon .%
\end{array}%
\end{equation*}
\end{definition}

\noindent Denote $\mathbb{F}^{\infty }=\mathbb{F}_{l}^{\infty }\cap \mathbb{F%
}_{r}^{\infty }$ the set of continuous non-anticipative functionals.

\begin{definition}[Boundedness-preserving functionals]
Define $\mathbb{B}$ as the set of non-anticipative functionals $F$ such that
for any compact $K\subset U$ and any $R>0$, there exists $C_{K,R}>0$ such
that:%
\begin{equation*}
\begin{array}{c}
\forall t\leq T,\text{ }\forall \left( x,v\right) \in D\left( \left[ 0,t%
\right] ,K\right) \times \mathcal{S}_{t}, \\
\sup\limits_{s\in \left[ 0,t\right] }\left\vert v\left( s\right) \right\vert
<R\Rightarrow \left\vert F_{t}\left( x,v\right) \right\vert <C_{K,R}.%
\end{array}%
\end{equation*}
\end{definition}

Next we introduce two notions of pathwise derivatives for a non-anticipative
functional $F=\left( F_{t}\right) _{t\in \left[ 0,T\right] }$: the
horizontal derivative and vertical derivative, respectively.

\begin{definition}[Horizontal derivative]
The horizontal derivative at $\left( x,v\right) \in \mathcal{U}_{t}\times
\mathcal{S}_{t}$ of non-anticipative functional $F=\left( F_{t}\right)
_{t\in \left[ 0,T\right] }$ is defined as
\begin{equation}
\mathcal{D}_{t}F\left( x,v\right) =\lim\limits_{h\rightarrow 0+}\frac{%
F_{t+h}\left( x_{t,h},v_{t,h}\right) -F_{t}\left( x,v\right) }{h},
\label{2.2}
\end{equation}%
if the corresponding limit exists. If (\ref{2.2}) is defined for all $\left(
x,v\right) \in \Upsilon $ the map%
\begin{equation*}
\mathcal{D}_{t}F:\mathcal{U}_{t}\times \mathcal{S}_{t}\mapsto \mathbb{R}^{d},%
\text{ }\left( x,v\right) \rightarrow \mathcal{D}_{t}F\left( x,v\right)
\end{equation*}%
defines a non-anticipative functional $\mathcal{D}F=\left( \mathcal{D}%
_{t}F\right) _{t\in \left[ 0,T\right] },$ the horizontal derivative of $F.$
\end{definition}

\begin{definition}[Dupire derivative]
A non-anticipative functional $F=\left( F_{t}\right) _{t\in \left[ 0,T\right]
}$ is said to be vertically differentiable at $\left( x,v\right) \in D\left( %
\left[ 0,t\right] ,\mathbb{R}^{n}\right) \times D\left( \left[ 0,t\right] ,%
\mathbb{S}_{n}^{+}\right) $ if
\begin{equation*}
\begin{array}{lll}
\mathbb{R}^{n} & \rightarrow & \mathbb{R} \\
\text{ \ }e & \rightarrow & F_{t}\left( x_{t}^{e},v_{t}\right)%
\end{array}%
\end{equation*}%
is differentiable at $0.$ Its gradient at $0$%
\begin{equation*}
\nabla _{x}F_{t}\left( x,v\right) =\left( \partial _{i}F_{t}\left(
x,v\right) ,i=1,\cdots ,d\right)
\end{equation*}%
where
\begin{equation}
\partial _{i}F_{t}\left( x,v\right) =\lim\limits_{h\rightarrow 0}\frac{%
F_{t}\left( x_{t}^{he_{i}},v\right) -F_{t}\left( x,v\right) }{h}  \label{2.3}
\end{equation}%
is called the vertical derivative of $F_{t}$ at $\left( x,v\right) .$ If (%
\ref{2.3}) is defined for all $\left( x,v\right) \in \Upsilon $, the
vertical derivative
\begin{equation*}
\nabla _{x}F:\mathcal{U}_{t}\times \mathcal{S}_{t}\mapsto \mathbb{R}^{n},%
\text{ }\left( x,v\right) \rightarrow \nabla _{x}F_{t}\left( x,v\right)
\end{equation*}%
defines a non-anticipative functional $\nabla _{x}F=\left( \nabla
_{x}F_{t}\right) _{t\in \left[ 0,T\right] }$ with value in $\mathbb{R}^{n}.$
\end{definition}

\begin{remark}
\textsl{If a vertically differentiable functional satisfies (\ref{2.1}), its
vertical derivative also satisfies (\ref{2.1}).}
\end{remark}

\begin{remark}
I\textsl{f }$F_{t}\left( x,v\right) =f\left( t,x\left( t\right) \right) $%
\textsl{\ with }$f\in C^{1,1}\left( \left[ 0,T\right) \times \mathbb{R}%
^{n}\right) $\textsl{\ then we retrieve the usual partial derivatives:}%
\begin{equation*}
\mathcal{D}_{t}F_{t}\left( x,v\right) =\partial _{t}f\left( t,x\left(
t\right) \right) ,\qquad \nabla _{x}F_{t}\left( x,v\right) =\nabla
_{x}f\left( t,x\left( t\right) \right) .
\end{equation*}
\end{remark}

\begin{remark}
\textsl{Note that if }$F$\textsl{\ is predictable with respect to the second
variable entails that for any }$t\in \left[ 0,T\right] ,$\textsl{\ }$%
F_{t}\left( x_{t},v_{t}^{e}\right) =F_{t}\left( x_{t},v_{t}\right) $\textsl{%
\ so an analogous notion of derivative with respect to }$v$\textsl{\ would
be identically zero under (\ref{2.1}).}
\end{remark}

\begin{definition}
Let $I\subset \left[ 0,T\right] $ be a subinterval of $\left[ 0,T\right] .$
Define $\mathbb{C}^{j,k}\left( I\right) $ as the set of non-anticipative
functionals $F=\left( F_{t}\right) _{t\in I}$ such that
\end{definition}

\begin{itemize}
\item $F$ is continuous at fixed times: $F_{t}:\mathcal{U}_{t}\times
\mathcal{S}_{t}\mapsto \mathbb{R}$ is continuous for the supremum norm.

\item $F$ admits $j$ horizontal derivative and $k$ vertical derivatives at
all $\left( x,v\right) \in \mathcal{U}_{t}\times \mathcal{S}_{t},$ $t\in I.$

\item $\mathcal{D}^{m}F$, $m\leq j,$ $\nabla _{x}^{n}F,$ $n\leq k$ are
continuous at fixed times $t\in I.$
\end{itemize}

\noindent We now introduce the following result (taken from Proposition 7 in
\cite{CF2})

\begin{proposition}[Functional Itô formula for a continuous semimartingale
on open set]
\label{pr1}Let $X$ be a continuous semimartingale valued in an open domain $%
U $ subset of $\mathbb{R}^{n}$ with quadratic variation process $\left[ X%
\right] ,$ and $A$ a continuous adapted process, on some filtered
probability space $\left( \tilde{\Omega},\mathcal{\tilde{F}},\mathcal{\tilde{%
F}}_{t},P\right) .$ Then for any non-anticipative functional $F\in \mathbb{C}%
^{1,2}\left( \left[ 0,T\right) \right) $ satisfying

(a) $F$ has a predictable dependence with respect to the second variable,
i.e. verifies (\ref{2.1}),

(b) $\nabla _{x}F,$ $\nabla _{x}^{2}F,$ $\mathcal{D}F\in \mathbb{B},$

(c) $F\in \mathbb{F}_{l}^{\infty },$

(d) $\nabla _{x}F,$ $\nabla _{x}^{2}F\in \mathbb{F}_{r}^{\infty },$

\noindent we have
\begin{eqnarray*}
F_{s}\left( X_{s},A_{s}\right) -F_{t}\left( X_{t},A_{t}\right)
&=&\int_{t}^{s}\mathcal{D}_{u}F\left( X_{u},A_{u}\right) \mbox{\rm d}u+\frac{%
1}{2}\int_{t}^{s}\text{Tr}\left[ ^{t}\nabla _{x}^{2}F_{u}\left(
X_{u},A_{u}\right) \mbox{\rm d}\left[ X\right] \left( u\right) \right] \\
&&+\int_{t}^{s}\nabla _{x}F_{u}\left( X_{u},A_{u}\right) \mbox{\rm d}X\left(
u\right) ,\text{ }P\text{-a.s.,}
\end{eqnarray*}%
where last term is the Itô's stochastic integral with respect to the $X.$
\end{proposition}

Now consider the continuous semimartingale $X$ driven as in (\ref{1.1}) and
denote
\begin{equation*}
\left[ X\right] \left( s\right) =\int_{t}^{s}\sigma \left( X_{r}\right)
\sigma ^{\ast }\left( X_{r}\right) \mbox{\rm d}r=\int_{t}^{s}A\left(
r\right) \mbox{\rm d}r
\end{equation*}%
its quadratic variation process. From now on, let $A\left( r\right) =\sigma
\left( X_{r}\right) \sigma ^{\ast }\left( X_{r}\right) $, which takes values
in the set $\mathbb{S}_{n}^{+}$ of symmetric positive $n\times n$ matrices,
has cadlag paths. We define the pair $\left( X_{r},A_{r}\right) =\left(
X_{r},\sigma \left( X_{r}\right) \sigma ^{\ast }\left( X_{r}\right) \right) $%
. Note that $A$ needs not to be a semimartingale. In particular, $F\left(
X_{u},A_{u}\right) =F\left( X_{u},A_{u-}\right) $ where $A_{u-}$ denotes the
path defined on $\left[ 0,t\right] $ by
\begin{equation*}
A_{u-}\left( r\right) =A\left( r\right) ,\text{ }r\in \left[ 0,u\right)
,\qquad A_{u-}\left( u\right) =A\left( u-\right) .
\end{equation*}%
For our purpose, let $U$ be an open subset of $\mathbb{R}^{n}$ with nonempty
boundary $\partial U$ and closure $\overline{U}.$ We introduce some useful
notations which are mainly taken and amended from \cite{P}.
\begin{eqnarray*}
\Lambda _{U_{t}} &\triangleq &\left\{ \gamma _{t}\in D\left( \left[ 0,t%
\right] ,\mathbb{R}^{n}\right) :\gamma \left( s\right) \in U,\text{ }s\in %
\left[ 0,t\right] \right\} , \\
\Lambda _{U} &\triangleq &\cup _{t\in \left[ 0,+\infty \right) }\Lambda
_{U_{t}}, \\
\Lambda _{\overline{U}_{t}} &\triangleq &\left\{ \gamma _{t}\in D\left( %
\left[ 0,t\right] ,\mathbb{R}^{n}\right) :\gamma \left( s\right) \in U,s\in %
\left[ 0,t\right) ,\gamma \left( t\right) \in \overline{U}\right\} , \\
\Lambda _{\overline{U}} &\triangleq &\cup _{t\in \left[ 0,+\infty \right)
}\Lambda _{\overline{U}_{t}}, \\
\Lambda _{\partial U_{t}} &\triangleq &\left\{ \gamma _{t}\in D\left( \left[
0,t\right] ,\mathbb{R}^{n}\right) :\gamma \left( s\right) \in U,s\in \left[
0,t\right) ,\gamma \left( t\right) \in \partial U\right\} , \\
\Lambda _{\partial U} &\triangleq &\cup _{t\in \left[ 0,+\infty \right)
}\Lambda _{\partial U_{t}}, \\
\Lambda &\triangleq &\cup _{t\in \left[ 0,+\infty \right) }D\left( \left[ 0,t%
\right] ,\mathbb{R}^{n}\right) .
\end{eqnarray*}

\noindent We now make the following assumptions for SDE (\ref{1.1}).

\begin{enumerate}
\item[\textbf{(A1)}] $\mu $ and $\sigma $ are two functionals satisfying%
\begin{equation*}
\left\{
\begin{array}{lllll}
\mu & : & \Lambda & \rightarrow & \mathbb{R}^{n}, \\
\sigma & : & \Lambda & \rightarrow & \mathbb{R}^{n\times d}.%
\end{array}%
\right.
\end{equation*}
\end{enumerate}

\begin{enumerate}
\item[\textbf{(A2)}] $\mu $ and $\sigma $ are Lipschitz in $x\in \Lambda ,$
i.e., there is some constant $C_{1}>0$ such that
\begin{equation*}
\left\vert \mu \left( x_{t}\right) -\mu \left( x_{t}^{\prime }\right)
\right\vert +\left\vert \sigma \left( x_{t}\right) -\sigma \left(
x_{t}^{\prime }\right) \right\vert \leq C_{1}\left\Vert x_{t}-x_{t}^{\prime
}\right\Vert ,\text{ for }x_{t},\text{ }x_{t}^{\prime }\in \Lambda ,
\end{equation*}%
and satisfy linear growth, i.e., there is some constant $C_{2}>0$ such that
\begin{equation*}
\left\vert \mu \left( x_{t}\right) \right\vert +\left\vert \sigma \left(
x_{t}\right) \right\vert \leq C_{2}\left( 1+\left\Vert x_{t}\right\Vert
\right) ,\text{ for }x_{t}\in \Lambda .
\end{equation*}
\end{enumerate}

\begin{remark}
\textsl{Note that }$\mu $\textsl{\ and }$\sigma $\textsl{\ are not necessary
non-anticipative functionals.}
\end{remark}

\begin{lemma}
Assume that the assumptions (A1)-(A2) hold. Then there exists a unique
strong solution to (\ref{1.1}).
\end{lemma}

\noindent The proof can be found in \cite{RY} Chapter IX, Theorem 2.1.

\noindent In this paper, we also use the following notations mainly taken
from \cite{CPF}.

\begin{itemize}
\item Let $C^{2,1}\left( A\right) $ denote all twice differentiable
functions on $A$, with bounded Lipschitz second order derivatives, where $A$
is an open subset of $\mathbb{R}^{n}.$

\item Let $S\subset \mathbb{R}^{n}$ be a nonempty set. We denote by $d_{S}$
the Euclidean distance function from $S,$ that is,
\begin{equation*}
d_{S}\left( x\right) =\inf\limits_{y\in S}\left\vert x-y\right\vert ,\qquad
\forall x\in \mathbb{R}^{n}.
\end{equation*}%
If $S$ is closed, then the above infimum is a minimum, which is attained on
a set that will be called the \textit{projection} of $x\in \mathbb{R}^{n}$
onto $S,$ i.e.
\begin{equation*}
\text{Proj}_{S}\left( x\right) =\left\{ y\in S:\left\vert x-y\right\vert
=d_{S}\left( x\right) \right\} ,\qquad \forall x\in \mathbb{R}^{n}.
\end{equation*}%
For every $t\in \left[ 0,+\infty \right) ,$ $x\in S,$ the hitting time of $%
\partial S$ is the random variable defined by
\begin{equation*}
\tau _{S}\left( \chi _{t}\right) =\inf \left\{ s\geq t:X^{\chi _{t}}\left(
s\right) \in \partial S\right\} ,\qquad \forall \chi _{t}\in \Lambda .
\end{equation*}

\item Let $\mathcal{K}$ be a closed subset of $\mathbb{R}^{n}$ with nonempty
interior $\overset{\circ }{\mathcal{K}}$ and boundary $\partial \mathcal{K}.$
We introduce the oriented distance function from $\partial \mathcal{K},$
i.e., the function
\begin{equation*}
b_{\mathcal{K}}\left( x\right) =\left\{
\begin{array}{ll}
d_{\partial \mathcal{K}}\left( x\right) & \text{if }x\in \overset{\circ }{%
\mathcal{K}}, \\
0 & \text{if }x\in \partial \mathcal{K}, \\
-d_{\partial \mathcal{K}}\left( x\right) & \text{if }x\in \mathcal{K}^{c},%
\end{array}%
\right.
\end{equation*}%
where $\mathcal{K}^{c}$ is the complimentary of $\mathcal{K}.$ In what
follows, we will use the following sets, defined for any $\varepsilon >0$%
\begin{eqnarray*}
\mathcal{N}_{\varepsilon } &\triangleq &\left\{ x\in \mathbb{R}%
^{n}:\left\vert b_{\mathcal{K}}\left( x\right) \right\vert <\varepsilon
\right\} , \\
\mathcal{K}_{\varepsilon } &\triangleq &\mathcal{K\cap N}_{\varepsilon }, \\
\overset{\circ }{\mathcal{K}_{\varepsilon }} &\triangleq &\overset{\circ }{%
\mathcal{K}}\mathcal{\cap N}_{\varepsilon }.
\end{eqnarray*}
\end{itemize}

\section{Main Result}

In this section, we study the viability properties of a compact $C^{2,1}$%
-smooth domain $\mathcal{K}$ defined as following with respect to the flow $%
X^{\chi _{t}}\left( \cdot \right) $ associated with the SDE (\ref{1.1})
(with Lipschitz continuous coefficients $\mu $ and $\sigma $ satisfying (A1)
and (A2)).

Keep in mind that $\mathcal{K}$ is a compact domain of class $C^{2,1}.$ From
Theorem 5.6 in \cite{DZ}, we know that
\begin{equation}
\mathcal{K}\text{ compact domain of class }C^{2,1}\Leftrightarrow \exists
\varepsilon _{0}>0\text{, such that }b_{\mathcal{K}}\in C^{2,1}\left(
\mathcal{N}_{\varepsilon _{0}}\right) .  \label{3.1}
\end{equation}%
A useful consequence of the above property is that
\begin{equation*}
\forall x\in \mathcal{K}_{\varepsilon _{0}},\text{ }\left\{
\begin{array}{l}
\exists !\text{ }\bar{x}\in \partial \mathcal{K}:b_{\mathcal{K}}\left(
x\right) =\left\vert x-\bar{x}\right\vert , \\
D_{x}b_{\mathcal{K}}\left( x\right) =D_{x}b_{\mathcal{K}}\left( \bar{x}%
\right) =-\vec{n}_{\mathcal{K}}\left( \bar{x}\right) ,%
\end{array}%
\right.
\end{equation*}%
where $\vec{n}_{\mathcal{K}}\left( \bar{x}\right) $ stands for the outward
unit normal to $\mathcal{K}$ at $\bar{x}.$ See, e.g. \cite{DZ}. Besides, if $%
\mathcal{K}$ is a compact domain of class $C^{2,1},$ one can see that there
is a sequence $\left\{ Q_{i}\right\} $ of compact domains of class $C^{2,1}$
such that
\begin{equation}
Q_{i}\subset \overset{\circ }{Q}_{i+1}\text{ and }\bigcup\limits_{i=1}^{%
\infty }Q_{i}=\overset{\circ }{\mathcal{K}}.  \label{3.2}
\end{equation}%
Indeed, it suffices to take that, for all $i$ large enough,
\begin{equation*}
Q_{i}=\left\{ x\in \mathcal{K}:b_{\mathcal{K}}\left( x\right) \geq \frac{1}{i%
}\right\} .
\end{equation*}%
Now suppose that
\begin{equation}
0\notin co\left\{ \nabla _{x}b_{\mathcal{K}}\left( x\right) \right\} ,\qquad
\forall x\in \partial \mathcal{K}.  \label{3.3}
\end{equation}%
Then, Clarke's tangent cone to $\mathcal{K}$ at every $x\in \mathcal{K}$ has
nonempty interior. For this reason, $\mathcal{K}$ coincides with the closure
of $\overset{\circ }{\mathcal{K}}$ (for more details see \cite{AF, CPF}).
Finally, we note that the existence of a sequence $\left\{ Q_{i}\right\} $
of compact domain of class $C^{2,1}$ satisfying (\ref{3.2}) is also
guaranteed when $\mathcal{K}$ is compact set with above properties (\ref{3.3}%
).

We now give our main result as follows:

\begin{theorem}
\label{th1} Suppose that assumptions (H1), (H2), and (\ref{3.3}) hold. Then
the following three statements are equivalent:

(i)\quad $\mathcal{K}$ enjoys the viability with respect to (\ref{1.1});

(ii) For $\forall \mathcal{X}_{t}\in \Lambda _{\partial \mathcal{K}_{t}},$
we have\quad
\begin{equation}
\left\{
\begin{array}{rcl}
\mathcal{L}_{\left( \mathcal{X}_{t},A_{t}\right) }\left( -b_{\mathcal{K}%
}\right) \left( \mathcal{X}_{t}\left( t\right) \right) & \leq & 0, \\
\sigma ^{\ast }\left( \mathcal{X}_{t}\right) \nabla _{x}b_{\mathcal{K}%
}\left( \mathcal{X}_{t}\left( t\right) \right) & = & 0,%
\end{array}%
\right.  \label{3.4}
\end{equation}%
where
\begin{equation*}
\mathcal{L}_{\left( \chi _{t,A_{t}}\right) }\varphi _{t}=\mathcal{D}%
_{t}\varphi +\left\langle \nabla _{x}\varphi _{t},\mu \right\rangle +\frac{1%
}{2}\text{Tr}\left( ^{t}\nabla _{x}^{2}\varphi _{t}A\left( t\right) \right) ,
\end{equation*}%
and $A\left( t\right) =\sigma \left( \mathcal{X}_{t}\right) \sigma ^{\ast
}\left( \mathcal{X}_{t}\right) ,$ for $\varphi \in \mathbb{C}^{2,1}\left( %
\left[ 0,T\right) \right) ;$

(iii)\quad $\overset{\circ }{\mathcal{K}}$ enjoys the viability with respect
to (\ref{1.1}).
\end{theorem}

\noindent For the proof of this theorem we need two auxiliary lemmata.

\begin{lemma}
\label{le1} Assume that (H1) and (H2) hold. Then for any compact domain $%
\mathbb{D}$ of $\mathbb{R}^{n}$ and for any $\chi _{t}\in \Lambda _{\mathbb{D%
}_{t}},$ the solution $X^{\chi _{t}}$ to (\ref{1.1}) satisfies
\begin{equation}
X^{\chi _{t}}\left( s\right) \in \mathbb{D},\text{ }s\in \left[ t,+\infty
\right) ,\text{ a.s..}  \label{3.5}
\end{equation}%
Then, for any $\mathbb{C}^{2,1}\left( \left[ 0,T\right) \right) $
non-anticipative functional $\varphi $ satisfying (a)-(d) in Proposition \ref%
{pr1} with a left frozen maximum (see Lemma 6 in \cite{P}) on $\Upsilon ^{%
\mathbb{D}\times \mathbb{S}_{n}^{+}}$ at $\left( \chi _{t},A_{t}\right) ,$
we have
\begin{equation*}
\left\{
\begin{array}{l}
\mathcal{D}_{t}\varphi \left( \chi _{t},A_{t}\right) +\left\langle \nabla
_{x}\varphi _{t}\left( \chi _{t},A_{t}\right) ,\mu \left( \chi _{t}\right)
\right\rangle +\frac{1}{2}\text{Tr}\left( ^{t}\nabla _{x}^{2}\varphi
_{t}\left( \chi _{t},A_{t}\right) A\left( t\right) \right) \leq 0, \\
\sigma ^{\ast }\left( \chi _{t}\right) \nabla _{x}\varphi _{t}\left( \chi
_{t},A_{t}\right) =0.%
\end{array}%
\right.
\end{equation*}
\end{lemma}

\proof%
Suppose that $\varphi $ is the non-anticipative test functional with a left
frozen maximum on $\Upsilon ^{\mathbb{D}\times \mathbb{S}_{n}^{+}}$ at $%
\left( \chi _{t},A_{t}\right) $ with $t\in \left[ 0,T\right) .$ More
precisely, for any given but fixed $\left( y_{s},z_{s}\right) \in \Upsilon ^{%
\mathbb{D}\times \mathbb{S}_{n}^{+}},$ $s\in \left[ t,T\right) ,$ we have
\begin{equation}
\varphi _{s}\left( \chi _{t}\oplus y_{s},A_{t}\oplus z_{s}\right) \leq
\varphi _{t}\left( \chi _{t},A_{t}\right) .  \label{3.6}
\end{equation}%
Immediately, from (\ref{3.5}) and (\ref{3.6}), we have
\begin{equation*}
\forall s\geq t,\text{ }\varphi _{s}\left( X_{s}^{\chi _{t}},A_{s}\right)
\leq \varphi _{t}\left( \chi _{t},A_{t}\right) ,\qquad P\text{-a.s.}.
\end{equation*}%
From standard arguments, the fact that this inequality holds for any $s\geq
t $, implies that, by virtue of Proposition \ref{pr1},
\begin{equation}
\mathcal{D}_{t}\varphi \left( \chi _{t},A_{t}\right) +\left\langle \nabla
_{x}\varphi _{t}\left( \chi _{t},A_{t}\right) ,\mu \left( \chi _{t}\right)
\right\rangle +\frac{1}{2}\text{Tr}\left( ^{t}\nabla _{x}^{2}\varphi
_{t}\left( \chi _{t},A_{t}\right) A\left( t\right) \right) \leq 0.
\label{3.7}
\end{equation}%
Let $\beta :\mathbb{R}\mathbf{\rightarrow }\mathbb{R}$ be an increasing
function such that $\beta ^{^{\prime }}\left( \varphi _{t}\left( \chi
_{t},A_{t}\right) \right) =1$ and $\beta ^{^{\prime \prime }}\left( \varphi
_{t}\left( \chi _{t},A_{t}\right) \right) =\lambda $ where $\lambda >0$ is
arbitrary$.$ It is easy to check that $\beta \circ \varphi _{t}$ also
attains a left frozen maximum on $\Upsilon ^{\mathbb{D}\times \mathbb{S}%
_{n}^{+}}$ at $\left( \chi _{t},A_{t}\right) .$ Hence, from (\ref{3.7}), we
have at $\left( \chi _{t},A_{t}\right) $ that%
\begin{equation}
\mathcal{D}_{t}\varphi +\left\langle \mu ,\nabla _{x}\varphi
_{t}\right\rangle +\frac{1}{2}\text{Tr}\left( \nabla _{x}^{2}\varphi
_{t}\sigma ^{\ast }\sigma \right) +\frac{\lambda }{2}\text{Tr}\left( \sigma
^{\ast }\sigma \nabla _{x}\varphi _{t}\nabla _{x}^{\ast }\varphi _{t}\right)
\leq 0,  \label{3.8}
\end{equation}%
where we have omitted the dependence in $\left( \chi _{t},A_{t}\right) $ for
simplicity. Obviously,
\begin{equation*}
\lambda \text{Tr}\left( \sigma ^{\ast }\sigma \left( \chi _{t}\right) \nabla
_{x}\varphi _{t}\left( \chi _{t},A_{t}\right) D_{x}^{\ast }\varphi
_{t}\left( \chi _{t},A_{t}\right) \right) =\lambda \left\vert \sigma ^{\ast
}\left( \chi _{t}\right) \nabla _{x}\varphi _{t}\left( \chi
_{t},A_{t}\right) \right\vert ^{2}.
\end{equation*}%
Noting that (\ref{3.8}) is still bounded as $\lambda \rightarrow +\infty ,$
we imply that
\begin{equation}
\sigma ^{\ast }\left( \chi _{t}\right) \nabla _{x}\varphi _{t}\left( \chi
_{t},A_{t}\right) =0.  \label{3.9}
\end{equation}%
The proof is complete.
\endproof%

\begin{remark}
\label{re1}\textsl{We do not impose the assumption (compact }$C^{2,1}$%
\textsl{-smooth domain) on }$\mathbb{D}$\textsl{\ in Lemma \ref{le1}}$.$
\end{remark}

As a consequence of Lemma \ref{le1}, we have the following:

\begin{corollary}
\label{co1}Suppose the same assumptions in Lemma \ref{le1} but take $%
\mathcal{K}$ instead of $\mathbb{D}$. Set
\begin{equation*}
\varphi \left( \overline{\chi }_{t},A_{t}\right) =-b_{\mathcal{K}}\left(
\overline{\chi }_{t}\left( t\right) \right) ,\text{ }\overline{\chi }_{t}\in
\Lambda _{\mathcal{K}_{t}},\text{ }t\in \left[ 0,+\infty \right) .
\end{equation*}%
Then we have, for $\forall \overline{\chi }_{t}\in \Lambda _{\partial
\mathcal{K}_{t}}$%
\begin{equation}
\left\{
\begin{array}{rcl}
\mathcal{L}_{\left( \overline{\chi }_{t},A_{t}\right) }\left( -b_{\mathcal{K}%
}\right) \left( \overline{\chi }_{t}\left( t\right) \right) & \leq & 0, \\
\sigma ^{\ast }\left( \overline{\chi }_{t}\right) \nabla _{x}b_{\mathcal{K}%
}\left( \overline{\chi }_{t}\left( t\right) \right) & = & 0.%
\end{array}%
\right.  \label{3.10}
\end{equation}
\end{corollary}

\proof%
Let us check that $\varphi \left( x_{t},A_{t}\right) =-b_{\mathcal{K}}\left(
x_{t}\left( t\right) \right) $ satisfies (b), (c), (d) in Proportion \ref%
{pr1}. First, (b) can be verified by assumption on $\mathcal{K}$ which is a
compact domain of class $C^{2,1}.$ For (c), let $x\in \Lambda _{\mathcal{K}%
_{t}},$ $x^{\prime }\in \Lambda _{\mathcal{K}_{t-h}}.$ Then, we have%
\begin{equation*}
\left\vert b_{\mathcal{K}}\left( x_{t}\left( t\right) \right) -b_{\mathcal{K}%
}\left( x_{t-h}^{\prime }\left( t-h\right) \right) \right\vert \leq
\left\vert x_{t}\left( t\right) -x_{t-h}^{\prime }\left( t-h\right)
\right\vert \leq d_{\infty }\left( x_{t},x_{t-h}^{\prime }\right) .
\end{equation*}%
Similarly, let $x\in \Lambda _{\mathcal{K}_{t}},$ $x^{\prime }\in \Lambda _{%
\mathcal{K}_{t+h}}.$ We also have%
\begin{equation*}
\left\vert b_{\mathcal{K}}\left( x_{t}\left( t\right) \right) -b_{\mathcal{K}%
}\left( x_{t+h}^{\prime }\left( t+h\right) \right) \right\vert \leq
\left\vert x_{t}\left( t\right) -x_{t+h}^{\prime }\left( t+h\right)
\right\vert \leq d_{\infty }\left( x_{t},x_{t+h}^{\prime }\right) ,
\end{equation*}%
since the Lipschitz constant of $b_{\mathcal{K}}$ is one. It remains to
check (d). In fact, by (\ref{3.1}) and the definition of $C^{2,1}\left(
A\right) ,$ $\nabla _{x}b_{\mathcal{K}},$ $\nabla _{x}^{2}b_{\mathcal{K}}$
can be also proved to be locally Lipschitz for the metric $d_{\infty }$.
\endproof%
%
%
%
%
%
%
%
%
%
%

Now let us turn back to compact $C^{2,1}$-smooth domain $\mathcal{K}.$ The
following approach is mainly borrowed from \cite{CPF}. It is not restrictive
to assume that $\varepsilon _{0}>0$ is such that, there exists a function $%
g\in C^{2,1}\left( \mathbb{R}^{n}\right) $ satisfying%
\begin{equation*}
\left\{
\begin{array}{ll}
0\leq g\leq 1, & \text{on }\mathcal{K}, \\
0<g, & \text{on }\mathcal{K}\backslash \mathcal{K}_{\varepsilon _{0}}, \\
g\equiv b_{\mathcal{K}}, & \text{on }\mathcal{K}_{\varepsilon _{0}}.%
\end{array}%
\right.
\end{equation*}%
Now define%
\begin{equation*}
\Psi \left( x_{t},A_{t}\right) =-\log \left( g\left( x_{t}\left( t\right)
\right) \right) ,\qquad \forall x_{t}\in \Lambda _{\overset{\circ }{\mathcal{%
K}}_{t}}.
\end{equation*}

\begin{lemma}
\label{le2}Assume that (H1), (H2), and (\ref{3.3}) hold. Then, there exists
a positive constant $M>0,$ such that%
\begin{equation}
\mathcal{L}_{\left( \chi _{t},A_{t}\right) }\Psi \left( \chi
_{t},A_{t}\right) \leq M,\text{ }\forall \chi _{t}\in \Lambda _{\overset{%
\circ }{\mathcal{K}}_{t}},\text{ }\forall t\in \mathbb{R}^{+},  \label{3.11}
\end{equation}%
where
\begin{equation*}
\mathcal{L}_{\left( \chi _{t},A_{t}\right) }\Psi \left( \chi
_{t},A_{t}\right) =\frac{1}{g^{2}\left( \chi _{t}\left( t\right) \right) }%
\left\vert \sigma ^{\ast }\left( \chi _{t}\right) \nabla _{x}g\left( \chi
_{t}\left( t\right) \right) \right\vert ^{2}-\frac{1}{g\left( \chi \left(
t\right) \right) }\mathcal{L}_{\left( \chi _{t},A_{t}\right) }g\left( \chi
_{t}\left( t\right) \right) .
\end{equation*}
\end{lemma}

\proof%
To begin with, we first give some estimations which will be very useful in
the sequel. We first notice that
\begin{equation*}
\frac{\partial }{\partial t}\Psi \left( x_{t},A_{t}\right) \equiv 0\text{
and }\Psi \left( x_{t},A_{t}\right) =\lim\limits_{\Lambda _{\overset{\circ }{%
\mathcal{K}}_{t}}\ni y_{t}\overset{d_{\infty }}{\rightarrow }x_{t}}\Psi
\left( y_{t},A_{t}\right) =+\infty ,\text{ }\forall x_{t}\in \Lambda
_{\partial \mathcal{K}_{t}}\text{.}
\end{equation*}%
Then, after a simple calculation, we have
\begin{equation*}
\mathcal{L}_{\left( \chi _{t},A_{t}\right) }\Psi \left( \chi
_{t},A_{t}\right) =\frac{1}{g^{2}\left( \chi _{t}\left( t\right) \right) }%
\left\vert \sigma ^{\ast }\left( \chi _{t}\right) D_{x}g\left( \chi
_{t}\left( t\right) \right) \right\vert ^{2}-\frac{1}{g\left( \chi
_{t}\left( t\right) \right) }\mathcal{L}_{\left( \chi _{t},A_{t}\right)
}g\left( \chi _{t}\left( t\right) \right) ,\text{ }\chi _{t}\in \Lambda _{%
\mathcal{K}_{t}}.
\end{equation*}%
We claim that, for any $\chi _{t}\in \Lambda _{\overset{\circ }{\mathcal{K}}%
_{t}},$
\begin{equation}
\frac{1}{g^{2}\left( \chi _{t}\left( t\right) \right) }\left\vert \sigma
^{\ast }\left( \chi _{t}\right) D_{x}g\left( \chi _{t}\left( t\right)
\right) \right\vert ^{2}-\frac{1}{g\left( \chi _{t}\left( t\right) \right) }%
\mathcal{L}_{\left( \chi _{t},A_{t}\right) }g\left( \chi _{t}\left( t\right)
\right) \leq M,\qquad \chi _{t}\left( t\right) \in \overset{\circ }{\mathcal{%
K}},  \label{3.12}
\end{equation}%
for some constant $M\geq 0.$ Indeed, the above estimate holds true$,$ when $%
\chi _{t}\left( t\right) \in \mathcal{N}_{\varepsilon _{0}}^{c},$ $b_{%
\mathcal{K}}\left( x\right) \geq \varepsilon _{0}$ since $g$ is strict
positive on $\mathcal{N}_{\varepsilon _{0}}^{c}$ and the assumption (H2),
where $\mathcal{N}_{\varepsilon _{0}}^{c}=\left\{ x\in \overset{\circ }{%
\mathcal{K}}:b_{\mathcal{K}}\left( x\right) \geq \varepsilon _{0}\right\} .$
Hence, we have to show (\ref{3.12}) holds for all $\chi _{t}\in \Lambda _{%
\overset{\circ }{\mathcal{K}}_{t}}$ satisfying $\chi _{t}\left( t\right) \in
\overset{\circ }{\mathcal{K}}\cap \mathcal{K}_{\varepsilon _{0}},$ i.e.,
\begin{equation*}
\frac{1}{b_{\mathcal{K}}^{2}\left( \chi _{t}\left( t\right) \right) }%
\left\vert \sigma ^{\ast }\left( \chi _{t}\right) D_{x}b_{\mathcal{K}}\left(
\chi _{t}\left( t\right) \right) \right\vert ^{2}-\frac{1}{b_{\mathcal{K}%
}\left( \chi _{t}\left( t\right) \right) }\mathcal{L}_{\left( \chi
_{t},A_{t}\right) }b_{\mathcal{K}}\left( \chi _{t}\left( t\right) \right)
\leq M,\qquad \forall \chi _{t}\left( t\right) \in \overset{\circ }{\mathcal{%
K}}\cap \mathcal{K}_{\varepsilon _{0}}.
\end{equation*}%
Given an $\chi _{t}\in \Lambda _{\overset{\circ }{\mathcal{K}}_{t}}$
satisfying $\chi _{t}\left( t\right) \in \overset{\circ }{\mathcal{K}}\cap
\mathcal{K}_{\varepsilon _{0}},$ let $\bar{x}\left( t\right) $ denote the
unique projection of $\chi _{t}\left( t\right) \in \overset{\circ }{\mathcal{%
K}}\cap \mathcal{K}_{\varepsilon _{0}}$ on the boundary of $\mathcal{K}$
since (\ref{3.1}). Then owning to (\ref{3.10}), we have
\begin{equation}
\sigma ^{\ast }\left( \chi _{t}^{\bar{x}\left( t\right) }\right) \nabla
_{x}b_{\mathcal{K}}\left( \bar{x}\left( t\right) \right) =0,  \label{3.13}
\end{equation}%
since $\chi _{t}^{\bar{x}\left( t\right) }\in \Lambda _{\partial \mathcal{K}%
_{t}}$. Therefore, we have
\begin{eqnarray*}
\left\vert \sigma ^{\ast }\left( \chi _{t}\right) \nabla _{x}b_{\mathcal{K}%
}\left( \chi _{t}\left( t\right) \right) \right\vert &=&\left\vert \left(
\sigma ^{\ast }\left( \chi _{t}\right) -\sigma ^{\ast }\left( \chi _{t}^{%
\bar{x}\left( t\right) }\right) \right) \nabla _{x}b_{\mathcal{K}}\left(
\chi _{t}\left( t\right) \right) +\sigma ^{\ast }\left( \chi _{t}^{\bar{x}%
\left( t\right) }\right) \nabla _{x}b_{\mathcal{K}}\left( \bar{x}\left(
t\right) \right) \right\vert \\
&=&\left\vert \left( \sigma ^{\ast }\left( \chi _{t}\right) -\sigma ^{\ast
}\left( \chi _{t}^{\bar{x}\left( t\right) }\right) \right) \nabla _{x}b_{%
\mathcal{K}}\left( \chi _{t}\left( t\right) \right) \right\vert \\
&\leq &\left\vert \sigma ^{\ast }\left( \chi _{t}\right) -\sigma ^{\ast
}\left( \chi _{t}^{\bar{x}\left( t\right) }\right) \right\vert \\
&\leq &C\left\Vert \chi _{t}-\chi _{t}^{\bar{x}\left( t\right) }\right\Vert
=C\left\vert \chi _{t}\left( t\right) -\bar{x}\left( t\right) \right\vert
=Cb_{\mathcal{K}}\left( \chi _{t}\left( t\right) \right) ,
\end{eqnarray*}%
where $C$ is a Lipschitz constant for $\sigma $ in the metric $\left\Vert
\cdot \right\Vert .$ Consequently,
\begin{equation}
\frac{1}{b_{\mathcal{K}}^{2}\left( \chi _{t}\left( t\right) \right) }%
\left\vert \sigma ^{\ast }\left( \chi _{t}\right) \nabla _{x}b_{\mathcal{K}%
}\left( \chi _{t}\left( t\right) \right) \right\vert ^{2}\leq \frac{\left(
Cb_{\mathcal{K}}\left( \chi _{t}\left( t\right) \right) \right) ^{2}}{b_{%
\mathcal{K}}^{2}\left( \chi _{t}\left( t\right) \right) }\leq C^{2}.
\label{3.14}
\end{equation}%
Also, noting (\ref{3.1}), we observe that,
\begin{equation*}
\mathcal{L}_{\left( \chi _{t},A_{t}\right) }b_{\mathcal{K}}\left( \chi
_{t}\left( t\right) \right) =\frac{1}{2}\text{Tr}\left( ^{t}\nabla
_{x}^{2}b_{\mathcal{K}}\left( \chi _{t}\left( t\right) \right) A\left(
t\right) \right) +\left\langle \mu \left( \chi _{t}\right) ,\nabla _{x}b_{%
\mathcal{K}}\left( \chi _{t}\left( t\right) \right) \right\rangle
\end{equation*}%
is Lipschitz continuous in $\overset{\circ }{\mathcal{K}}\cap \mathcal{K}%
_{\varepsilon _{0}}$ for the metric $\left\Vert \cdot \right\Vert $. Thus
\begin{equation}
\mathcal{L}_{\left( \chi _{t}^{\bar{x}\left( t\right) },A_{t}\right) }\left(
-b_{\mathcal{K}}\right) \left( \bar{x}\left( t\right) \right) \leq 0
\label{3.15}
\end{equation}%
yields that
\begin{eqnarray}
&&-\frac{1}{b_{\mathcal{K}}\left( \chi _{t}\left( t\right) \right) }\mathcal{%
L}_{\left( \chi _{t},A_{t}\right) }b_{\mathcal{K}}\left( \chi _{t}\left(
t\right) \right) \\
&=&\frac{-\mathcal{L}_{\left( \chi _{t},A_{t}\right) }b_{\mathcal{K}}\left(
\chi _{t}\left( t\right) \right) +\mathcal{L}_{\left( \chi _{t}^{\bar{x}%
\left( t\right) },A_{t}\right) }b_{\mathcal{K}}\left( \bar{x}\left( t\right)
\right) +\mathcal{L}_{\left( \chi _{t}^{\bar{x}\left( t\right)
},A_{t}\right) }\left( -b_{\mathcal{K}}\left( \bar{x}\left( t\right) \right)
\right) }{b_{\mathcal{K}}\left( \chi _{t}\left( t\right) \right) }  \notag \\
&\leq &\frac{1}{b_{\mathcal{K}}\left( \chi _{t}\left( t\right) \right) }%
\left\vert \mathcal{L}_{\left( \chi _{t}^{\bar{x}\left( t\right)
},A_{t}\right) }b_{\mathcal{K}}\left( \bar{x}\left( t\right) \right) -%
\mathcal{L}_{\left( \chi _{t},A_{t}\right) }b_{\mathcal{K}}\left( \chi
_{t}\left( t\right) \right) \right\vert  \notag \\
&\leq &C\frac{\left\Vert \chi _{t}-\chi _{t}^{\bar{x}\left( t\right)
}\right\Vert }{b_{\mathcal{K}}\left( \chi \left( t\right) \right) }=C\frac{%
\left\vert \chi _{t}\left( t\right) -\bar{x}\left( t\right) \right\vert }{b_{%
\mathcal{K}}\left( \chi \left( t\right) \right) }=C\frac{b_{\mathcal{K}%
}\left( \chi \left( t\right) \right) }{b_{\mathcal{K}}\left( \chi \left(
t\right) \right) }\leq C,  \label{3.16}
\end{eqnarray}%
for all $\chi _{t}\in \Lambda _{\overset{\circ }{\mathcal{K}}_{t}}$
satisfying $\chi _{t}\left( t\right) \in \overset{\circ }{\mathcal{K}}\cap
\mathcal{K}_{\varepsilon _{0}},$ where $C$ is a Lipschitz constant for $%
\mathcal{L}_{\left( \chi _{t},A_{t}\right) }b_{\mathcal{K}}.$ Consequently,
combining (\ref{3.14}) and (\ref{3.16}), we get desired result.
\endproof%

\begin{remark}
\label{re3}\textsl{Clearly, the estimations of (\ref{3.14}) and (\ref{3.16})
indicate that the initial condition }$\chi _{t}$\textsl{\ must be chosen in }%
$D\left( \left[ 0,t\right] ,\mathcal{K}\right) $\textsl{\ in Corollary \ref%
{co1}.}
\end{remark}

\noindent Now we are able to give the proof of Theorem \ref{th1}, similarly
as in \cite{CPF}.

\proof%
We first prove (i)$\Rightarrow $(ii). By Corollary \ref{co1}, the first
assertion holds.\textbf{\ }Next, we are going to show (ii)$\Rightarrow $%
(iii). Consider the stopping time $\tau _{Q_{i}}\left( \chi _{t}\right) $
where $\left\{ Q_{i}\right\} $ is the sequence of compact domain of class $%
C^{2,1}$ satisfying (\ref{3.2}). Applying the functional Itô formula
(Proposition \ref{pr1}) we get, for all $\chi _{t}\in \Lambda _{Q_{i,t}}$
and $0\leq t\leq s,$%
\begin{eqnarray}
&&\Psi \left( X_{s\wedge \tau _{Q_{i}}\left( \chi _{t}\right) }^{\chi
_{t}},A_{s\wedge \tau _{Q_{i}}\left( \chi _{t}\right) }\right) -\Psi \left(
\chi _{t},A_{t}\right)  \notag \\
&=&\int_{t}^{s\wedge \tau _{Q_{i,t}}\left( \chi _{t}\right) }\left( \mathcal{%
L}_{\left( X_{r}^{\chi _{t}},A_{r}\right) }\Psi \right) \left( X_{r}^{\chi
_{t}},A_{r}\right) \mbox{\rm d}r  \notag \\
&&+\int_{t}^{s\wedge \tau _{Q_{i}}\left( \chi _{t}\right) }\left\langle
\nabla _{x}\Psi \left( X_{r}^{\chi _{t}},A_{r}\right) ,\sigma \left(
X_{r}^{\chi _{t}}\right) \right\rangle \mbox{\rm d}W\left( r\right) .
\label{3.17}
\end{eqnarray}%
Then, taking expectation and noting (\ref{3.11}), we get
\begin{eqnarray*}
\mathbb{E}\left[ \Psi \left( X_{s\wedge \tau _{Q_{i}}\left( \chi _{t}\right)
}^{\chi _{t}},A_{s\wedge \tau _{Q_{i}}\left( \chi _{t}\right) }\right) %
\right] &=&\Psi \left( \chi _{t},A_{t}\right) +\mathbb{E}\left[
\int_{t}^{s\wedge \tau _{Q_{i}}\left( \chi _{t}\right) }\left( \mathcal{L}%
_{\left( X_{r}^{\chi _{t}},A_{r}\right) }\Psi \right) \left( X_{r}^{\chi
_{t}},A_{r}\right) \mbox{\rm d}r\right] \\
&\leq &\Psi \left( \chi _{t},A_{t}\right) +Ms.
\end{eqnarray*}%
By Fatou's Lemma, the above inequality yields that
\begin{equation*}
\mathbb{E}\left[ \Psi \left( X_{s\wedge \tau _{_{\mathcal{K}}}\left( \chi
_{t}\right) }^{\chi _{t}},A_{s\wedge \tau _{_{\mathcal{K}}}\left( \chi
_{t}\right) }\right) \right] \leq \Psi \left( \chi _{t},A_{t}\right) +Ms,%
\text{ }\forall s\geq t\geq 0,\qquad \forall \chi _{t}\in \Lambda _{\overset{%
\circ }{\mathcal{K}}_{t}}.
\end{equation*}%
Observing that the function in the right-hand above is finite on $\Lambda _{%
\overset{\circ }{\mathcal{K}}},$ we deduce that
\begin{equation*}
P\left( \tau _{_{\mathcal{K}}}\left( \chi _{t}\right) \leq s\right) =P\left(
\Psi \left( X_{s\wedge \tau _{_{\mathcal{K}}}\left( \chi _{t}\right) }^{\chi
_{t}},A_{s\wedge \tau _{_{\mathcal{K}}}\left( \chi _{t}\right) }\right)
=+\infty \right) =0,\text{ }s\geq t\geq 0,\qquad \forall \chi _{t}\in
\Lambda _{\overset{\circ }{\mathcal{K}}_{t}}.
\end{equation*}%
Take a sequence $t_{k}\uparrow +\infty $ and observe that
\begin{equation*}
0=P\left( \tau _{_{\mathcal{K}}}\left( \chi _{t}\right) \leq t_{k}\right)
\uparrow P\left( \tau _{_{\mathcal{K}}}\left( \chi _{t}\right) <+\infty
\right) ,\qquad \forall \chi _{t}\in \Lambda _{\overset{\circ }{\mathcal{K}}%
_{t}}.
\end{equation*}%
We complete the proof of (ii)$\Rightarrow $(iii).

Next we are going to prove (iii)$\Rightarrow $(i). Assume that $\overset{%
\circ }{\mathcal{K}}$ is viable and fix $\chi \in \Lambda _{\mathcal{K}_{t}}$%
. Let $\left\{ \chi _{t,k}\right\} $ be a sequence in $\Lambda _{\overset{%
\circ }{\mathcal{K}}_{t}}$ such that $\chi _{t,k}\overset{\left\Vert \cdot
\right\Vert }{\rightarrow }\chi _{t}$ since $\mathcal{K}$ coincides with the
closure of $\overset{\circ }{\mathcal{K}}.$ Then we have $X^{\chi
_{t,k}}\left( s\right) \in \mathcal{K}$, $P$-a.s. for all $s\geq t.$ Noting $%
X^{\chi _{t,k}}\left( s\right) \rightarrow X^{\chi _{t}}\left( s\right) ,$ $%
P $-a.s. for $s\geq t$, we deduce that
\begin{equation*}
X^{\chi _{t}}\left( s\right) \in \mathcal{K},\text{ }P\text{-a.s., for all }%
s\geq t.
\end{equation*}%
From the arbitrary point $\chi _{t}$ in $\Lambda _{\mathcal{K}_{t}}$, we get
the desired result.
\endproof%

\begin{remark}
\label{re4}\textsl{If we suppose that }$t=0,$\textsl{\ }$\chi \left(
0\right) =x,$\textsl{\ }$k\left( \gamma _{s}\right) =k\left( \gamma
_{s}\left( s\right) \right) $\textsl{\ where }$\gamma \in \Lambda ,$\textsl{%
\ }$x\in \mathcal{K},$\textsl{\ }$s\in \left[ 0,+\infty \right) ,$\textsl{\
while }$k=\mu ,$\textsl{\ }$\sigma ,$\textsl{\ respectively}$.$\textsl{\
Then we recover Theorem 3.2 in \cite{CPF}.}
\end{remark}

\noindent \textbf{Acknowledgments. }The author would like to thank two
anonymous referees, Professor Marc Quincampoix for their valuable comments,
which led to a much better version of this article.

\end{document}